\documentclass[12pt,a4paper]{article}
\usepackage[centertags]{amsmath}
\usepackage{amsfonts}

\usepackage{graphics}
\usepackage{amssymb}
\usepackage{amsthm}
\usepackage{newlfont}
\usepackage{epsfig}
\usepackage{amssymb}
\usepackage{amsmath}
\usepackage{amsthm}
\usepackage{graphicx}
\usepackage{makeidx}
\usepackage[mathscr]{euscript}
\theoremstyle{plain}
\newtheorem{thm}{Theorem}[section]

\newtheorem{lem}[thm]{Lemma}
\newtheorem{example}[thm]{Example}

\theoremstyle{definition}
\newtheorem{defn}{Definition}[section]

\theoremstyle{remark} \tolerance=10000 \hbadness=10000
\def \ni{\noindent}
	\author{Jeepamol J. Palathingal \footnote{Email : jeepamoljp@gmail.com}\\ Department of Mathematics\\
		PM Government College, Chalakudy-680722\\ \vspace{0.3cm} Kerala, India.\\
		Aparna Lakshmanan S.\footnote{E-mail :
			aparnaren@gmail.com, aparnals@cusat.ac.in}\\
		Department of Mathematics\\Cochin University of Science and Technology\\ \vspace{0.3cm} Cochin-682022, Kerala, India.\\
		Greg Markowsky \footnote{E-mail : greg.markowsky@monash.edu} \\
			School of Mathematics \\ Monash University \\Clayton, VIC Australia}
	
	\title{\textbf{The Gallai and anti-Gallai graphs of strongly regular graphs}}

	\date{}
	
\begin{document}
		
		\maketitle
		
		\begin{abstract}
			
			  In this paper, we show that if $G$ is strongly regular then the Gallai graph $\Gamma(G)$ and the anti-Gallai graph $\Delta(G)$ of $G$ are edge-regular. We also identify conditions under which the Gallai and anti-Gallai graphs are themselves strongly regular, as well as conditions under which they are 2-connected. We include also a number of concrete examples and a discussion of spectral properties of the Gallai and anti-Gallai graphs. \\

			\ni {\bf Keywords:} Gallai Graph, Anti-Gallai Graph, Strongly Regular Graph, Edge-regular Graph, Adjacency Spectrum\\

\ni{\bf AMS classification:} 05C50, 05C76.\\
			
			
		\end{abstract}
		\section{Introduction}

		Given a graph $G$, the {\it line graph} $L(G)$ of $G$ has the edges of
		$G$ as its vertices, with two vertices adjacent in
		$L(G)$ if the corresponding edges have a vertex in common in $G$. \if2 The adjacency spectrum of line graph of many classes of graphs has been studied in \cite{cv1},\cite{cv2} and \cite{IG}.\fi The line graph $L(G)$ can then be decomposed into the Gallai graph $\Gamma(G)$ and anti-Gallai graph $\Delta(G)$, as follows. The {\it Gallai graph} $\Gamma(G)$ of $G$ has the edges of $G$ as its vertices, with two vertices adjacent in
		$\Gamma(G)$ if their corresponding edges have a vertex in common in $G$ but do not lie on a common
		triangle in $G$. The {\it anti-Gallai graph} $\Delta(G)$ of $G$ again has
		the edges of $G$ as its vertices, with two vertices in $\Delta(G)$ adjacent if their corresponding edges lie on a common triangle in $ G$. We will refer to $G$ as the {\it underlying graph} of the graphs $L(G), \Gamma(G),$ and $\Delta(G)$. \\
		
		Clearly, $L(G)$ is the edge-disjoint union of $\Gamma(G)$ and $\Delta(G)$. Line graphs have been extensively studied in the literature, however in comparison the theory of Gallai and anti-Gallai graphs has attracted less attention since their introduction in \cite{gal}. Good overviews on the topic can be found in \cite{Le1} and \cite{Le2}, while more recent papers include \cite{Apa1,Apa,garg,Jee1,joos,agnes1,agnes2}. Furthermore, anti-Gallai graph have found application in linguistics, specifically in identifying polysemous words \cite{PHRC}. This provides motivation for examining the spectral properties of anti-Gallai graphs, particularly in relation to the spectra of their underlying graphs. The spectra of some classes of Gallai and anti-Gallai graphs have already been studied, in \cite{Jee3} and \cite{Jee2}.\\
		
		Our interest in this paper is in furthering the study of Gallai and anti-Gallai graphs, particularly in relating spectral and other properties of these graphs with those of their underlying ones. However, there is an issue which presents significant difficulties. Regularity properties of $G$ do not always translate to similar properties of $\Gamma(G)$ and $\Delta(G)$. In particular, $G$ being regular does not imply that $\Gamma(G)$ and $\Delta(G)$ are also regular. Contrast this with the situation for line graphs, where regularity of $G$ implies regularity of $L(G)$, a fact which has greatly facilitated the study of line graphs. However, a key observation in this paper is that $\Gamma(G)$ and $\Delta(G)$ of a strongly regular graph $G$ are regular, and in fact are even edge-regular. This has prompted us to focus on the Gallai and anti-Gallai graphs of strongly regular graphs, and this is the topic of this paper.\\
		
		We will now define strongly regular graphs. Let $G$ be a $k$-regular graph with $n$ vertices. The graph $G$ is said to be {\it strongly regular} \cite{Bap} with parameters $(n,k,\lambda,\mu)$ if the following conditions hold:\\
		
		\noindent (1) $G$ is neither complete nor empty;\\
		(2) any two adjacent vertices of $G$ have $\lambda$ common neighbours;\\
		(3) any two non-adjacent vertices of $G$ have $\mu$ common neighbours.\\

		We note that the cycles $C_4$ and $C_5$ are strongly regular, however all questions that we will consider reduce to trivialities in these cases, so we will henceforth assume $k \geq 3$ whenever $G$ is strongly regular. \\
		
		A weaker notion is that of an edge-regular graph. An {\it edge-regular graph} \cite{L} with parameters $(n,k,\lambda)$ is a graph on $n$ vertices which is regular of degree $k$ and such that any two adjacent vertices have exactly $\lambda$ common neighbours. It is evident that strongly regular graphs are edge-regular, but the converse does not hold. In fact, the class of strongly regular graphs is far more restricted than that of edge-regular ones, with a rich structure that has attracted numerous researchers. For excellent overviews on the topic, see \cite{Bro, cam,godroy}.  \\
		
		In the next section we give some necessary preliminaries on strongly regular graphs, and the subsequent section concerns the Gallai and anti-Gallai graphs of strongly regular graphs. All graph theoretic notations and terminology not defined here
		are standard, and can be found for instance in \cite{cv1} and \cite{Bal}.

\section{Preliminaries on strongly regular graphs}

Henceforth, we assume $G$ is a connected strongly regular graph with parameters $(n,k,\lambda,\mu)$. We will in many cases refer to subgraphs of $G$, and we must be careful since there are two different types of subgraphs which we will use. The first is an {\it induced subgraph}, which is formed from a subset of the vertices of $G$ together with all edges connecting vertices in that set. The more general definition of subgraphs allows for edges from $G$ to be absent even when both endpoints are present in the subgraph. To avoid confusion, we will always refer to induced sugbgraphs as such, and if we simply consider a subgraph (without the word "induced" present) then it is not assumed to be induced. \\

The following lemmas contain facts that are certainly known, but the proof is included for the benefit of readers unfamiliar with the structure of strongly regular graphs and standard methods of proof in the field.  \\\\

\begin{lem} \label{newlemma}
\begin{enumerate}
\item  $\mu \geq 1$.

\item If $\mu >1$ then every pair of non-adjacent vertices must belong to at least one cycle of length $4$.

\item Given an edge $e=uv$ and a vertex $w$, either there is a vertex $y$ adjacent to all of $w,u,v$ or else $w,u,v$ must belong to a cycle $C_n$ of length at most $5$.

\item Given a pair of edges $uv$ and $xy$, either there is a vertex $z$ adjacent to all of $u,v,x,y$, or else both edges lie on a cycle $C_n$ of length at most $6$.

\item If $\mu = 1$ then diamonds and $C_4$'s are forbidden as induced subgraphs in $G$ (a diamond refers to two triangles sharing a common edge, i.e $K_4$ with an edge removed \cite{Jee1}).

\item A strongly regular graph is $2$-connected (that is, the removal of any vertex does not disconnect the graph).

\item If $\lambda \leq 1$ then diamonds and $K_4$'s are forbidden in $G$.

\end{enumerate}
\end{lem}

\begin{proof}
We will employ the following notation, which is not entirely standard in the field. Given vertex $x$, let $N_j(x) = \{y:d(x,y)=j\}$. More standard would be to use $\Gamma$ instead of $N$, but we will reserve $\Gamma$ for the Gallai graph. Then, for any $x$, the vertices of $G$ can be partitioned into the sets $\{x\}, N_1(x),$ and $N_2(x)$.

Proving {\it 1} is trivial, since $G$ is connected and of diameter 2.

As for {\it 2}, if $y \in N_2(x)$ and $\mu > 1$ then there are distinct $u,v \in N_1(x) \cap N_1(y)$ and then $xuyv$ is a cycle of length 4.

For {\it 3}, since $G$ has diameter 2 there are paths of length at most 2 from $w$ to $u$ and from $w$ to $v$. If these paths do not share vertices other than $w$ then they, together with $uv$, form a cycle of length at most 5. It is also possible, however, that the paths are of the form $wxu$ and $wxv$ for some vertex $x$, but then $x$ is adjacent to $w,u,v$, and the conclusion holds.

{\it 4} follows similarly, since we can find paths of length at most 2 between $u$ and $x$ and between $v$ and $y$, and then these paths can then either be used to form a cycle of length at most 6 or else yield a point adjacent to all of $u, x, v,y$.

For {\it 5}, these graphs cannot exist as induced subgraphs since both contain vertices of distance 2 from each other but with two paths of length 2 connecting them.

For {\it 6}, if we remove a vertex $w$ from $G$, then we must show that the subgraph induced on $N_1(w) \cup N_2(w)$ is connected. Choose $x \in N_2(w)$, and let $y \in N_1(w) \cup N_2(w)$. Since the diameter of $G$ is 2, there must be a path of length at most 2 from $x$ to $y$, and this path can't pass through $w$ since $d(w,x) = 2$. The result follows from this.

For {\it 7}, note that $\lambda=1$ implies that any adjacent points can have at most one common neighbor. Both of the graphs indicated contain adjacent points with 2 common neighbors, and therefore can't exist as subgraphs of $G$.
\end{proof}

\begin{figure}
\begin{centering}
\includegraphics[width=7cm,height=6cm,keepaspectratio]{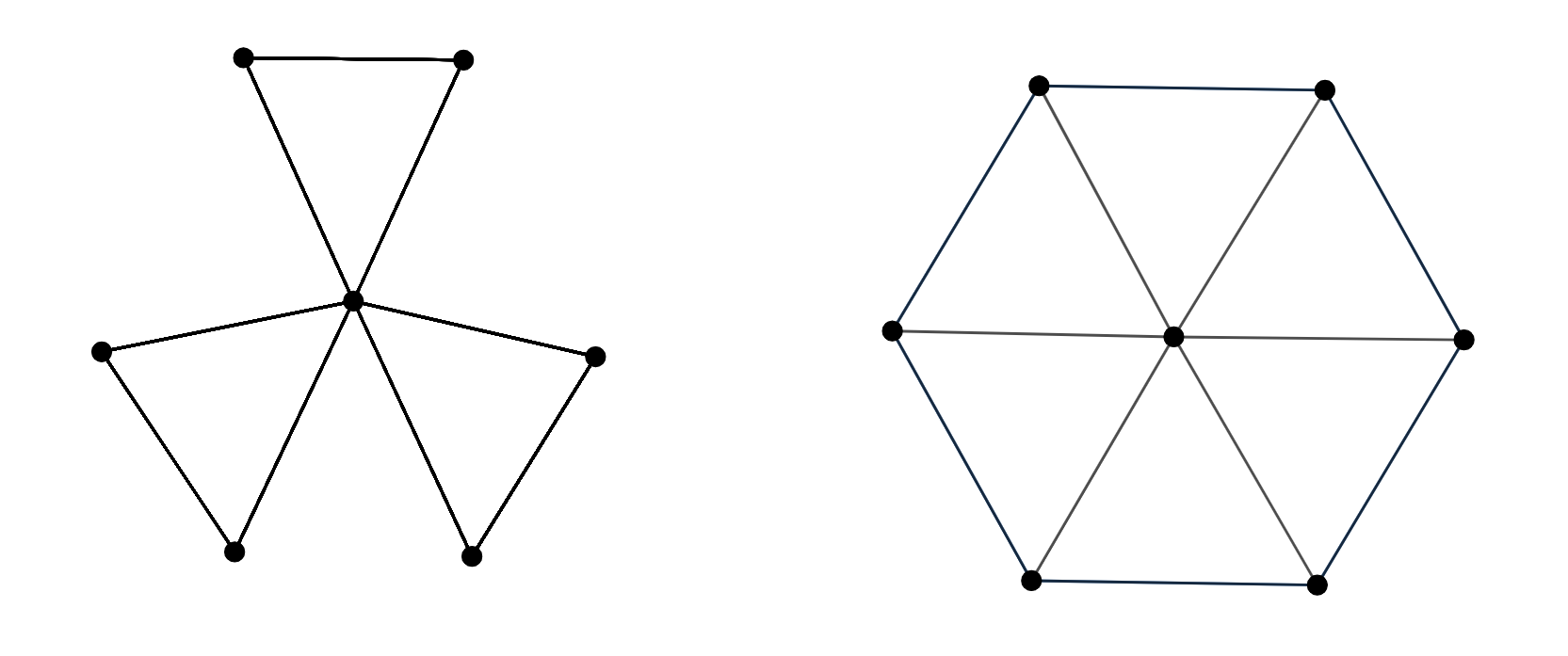}
\par\end{centering}
\caption{The fan graph $F_3$ and wheel graph $W_6$} \label{fanwheel}
\end{figure}

The {\it fan graph}, denoted by $F_n$, is $n$ copies of $K_3$ conjoined at a common vertex. Note that $F_n$ is a planar graph with $2n+1$ vertices and $3n$ edges \cite{Ka}. The following result gives a local characterization of strongly regular graphs with $\lambda=1$.

 \begin{lem}\label{11}
 	Let $G$ be a connected strongly regular graph with parameters $(n,k,1,\mu)$. Then every vertex together with its neighbours form a $k/2$-fan.
 \end{lem}
	\begin{proof}
	Consider an arbitrary vertex $u$ in $G$. Since $G$ is connected, there exists an edge incident to $u$ and by hypothesis that edge belongs to a $K_3$. $G \neq K_3$ and regular implies that there exists another edge incident to $u$ but disjoint from the first triangle. But by {\it 7} of the previous lemma,  diamonds are forbidden in $G$ and hence that edge also belongs to another $K_3$. Proceeding like this we find that $u$ together with its neighbours form a $k/2$-fan. \\
	\end{proof}

A {\it wheel graph} $W_n$ is formed by joining a single vertex to all vertices of the $n$-cycle $C_n$. The following result indicates a structural characteristic of strongly regular graphs with $\lambda \geq 2$.

\begin{lem}\label{12}
	Let $G$ be a connected strongly graph with parameters $(n,k,\lambda,\mu)$ with $\lambda \geq 2$. Then every vertex is a central vertex of at least one wheel.
\end{lem}

\begin{proof}
	Since $\lambda \geq 2$, given an edge $e$ the endpoints have at least 2 common neighbors, and $e$ together with these two paths of length 2 form a diamond. Consider $u\in V(G)$, and let $e_1$ be any edge incident to $u$. $e_1$ is then the central edge of a diamond. Let $e_2$ and $e_3$ be the other edges of the diamond which are incident to $u$. Then there exist diamonds with $e_2$ and $e_3$ as their central edges.  We proceed like this, and the process will terminate only when $u$ together with some set of edges form a wheel. It may happen that there exists another edge incident to $u$ which is not in the above wheel. In this case by the same argument we can find another wheel with $u$ as the central vertex.\\
\end{proof}	

\section{Gallai and anti-Gallai graphs of strongly regular graphs}

The following indicates the advantages enjoyed by requiring graphs to be edge-regular (or strongly regular) rather than just regular. We reiterate that this result is false for regular graphs which are not edge-regular.

\begin{thm}\label{15}
	If $G$ is an edge-regular graph with parameters $(n,k,\lambda)$, then $\Gamma(G)$ and $\Delta(G)$ are regular graphs.
\end{thm}
\begin{proof}
	Consider an arbitrary vertex $x$ of $\Gamma(G)$. Let $uv$ be the edge corresponding to this vertex in $G$. Then by the definition of $\Gamma(G)$, the degree of $x$ is
	$d_{\Gamma(G)}(x)= d_{G}(u)+d_{G}(v)- 2|N_1(u) \cap N_1(v)| -2$,
	where $d_{\Gamma(G)}(x)$ denotes the degree of $x$ in $\Gamma(G)$ and $d_{G}(u),d_{G}(v)$ denote the degrees of $u,v$ in $G$. Hence $\Gamma(G)$ is $2(k-\lambda-1)$-regular. \\

	Now consider an arbitrary vertex $x$ in $\Delta (G)$. Let $e=uv$ be the corresponding edge in $G$. Then the degree of $x$ in $\Delta(G)$ is $2|N_1(u) \cap N_1(v)| = 2\lambda$. Therefore, $\Delta(G)$ is $2\lambda$-regular.\\
\end{proof}

We are interested in even stronger statements. In particular, we are interested in conditions under which the Gallai and anti-Gallai graphs are edge-regular, or even strongly regular. We will consider each of these separately in the following two subsections.

\subsection{Gallai graphs}

As per Theorem \ref{15} we have that if $G$ is strongly regular then $\Gamma(G)$ is regular. In this section we prove that the Gallai graphs of some special classes of strongly regular graphs are edge-regular or strongly regular. We begin with a more general structure theorem. \\

\begin{thm}\label{14}
	Let $G$ be a connected graph. The Gallai graph $\Gamma(G)$ is disconnected if and only if there exists a partition of the edge set into $E_1, E_2, ... E_p$ where $p\geq2$, such that $e_i\in E_i$ and $e_j\in E_j$ are incident in $G$ implies $e_i$ and $e_j$ span a triangle in $G$.
\end{thm}
\begin{proof}
	Suppose that $\Gamma(G)$ is disconnected and let $\Gamma_1, \Gamma_2,...\Gamma_p$ with $p\geq 2$ be the components of $\Gamma(G)$. Let $E_i$=\{$e\in G$ : $e$ is an edge corresponding to a vertex $v$ in $\Gamma_i$\}, where $1\leq i \leq p$. Clearly, $E_i $ is a partition for  $E(G)$. Since the connectedness of $G$ implies the  connectedness of $L(G)$,  at least one edge  of $E_i$ is incident with some $e\in E_j$ for some $j \neq i$. But, $\Gamma_i$ and $\Gamma_j$ are different components of $\Gamma(G)$, and hence if $e_i\in E_i$ is incident with $e_j\in E_j$ then they must span a triangle in $G$.\\

	 For the converse assume that such a partition exists for $E(G)$. Then, for any distinct $i$ and $j$, the vertices corresponding to the edges in $E_i$ and $E_j$ are in different components in $\Gamma(G)$.\\
	  \end{proof}

We can say more in special cases. For instance, for $\lambda = 0$ we have the following.

	 \begin{thm}\label{16}
	 	Let $G$ be a connected strongly regular graph with parameters $(n,k,0,\mu)$. Then $\Gamma(G)$ is connected and edge-regular. Also it is strongly regular if and only if the following conditions hold:\\
	 	$1$. If $\mu=1$ then any two non-adjacent edges belong to a common $C_5$.\\
	 	$2$. If $\mu>1$ then any two non-adjacent edges belong to a common $C_4$.
	 \end{thm}
 \begin{proof}
 Since $\lambda=0$, $G$ is $K_3$-free and thus $\Gamma(G)\cong L(G)$. Since $G$ is connected and $k$-regular, $L(G)$ is connected and $2k-2$ regular, and thus so is $\Gamma(G)$.\\

 Consider two adjacent vertices $x$ and $y$ in $L(G)$, and let $e_1$ and $e_2$ be the corresponding edges in $G$. Then $e_1$ and $e_2$ must have a common vertex in $G$. The number of common vertices of $x$ and $y$ in $L(G)$ is same as the number of edges incident on the common vertex of both $e_1$ and $e_2$. Since the graph is $k$-regular it is equal to $k-2$. Hence $\Gamma(G)$ is an edge-regular graph with parameters $(\frac{nk}{2}, 2k-2, k-2)$.\\

 Now consider two non-adjacent vertices  $x$ and $y$ in $L(G)$. Let $f_1$ and $f_2$ be the corresponding edges in $G$. Then $f_1$ and $f_2$ have no  common vertex in $G$. The number of common vertices of $x$ and $y$ in $L(G)$ is same as the number of edges adjacent to both $f_1$ and $f_2$. \\

 If $\mu=1$, by Lemma \ref{newlemma} part {\it 4}, $f_1$ and $f_2$ belong to a cycle of length at most 6 (the other possibility, that one point is adjacent to all endpoints of $f_1$ and $f_2$, is not possible since $\lambda = 0$). The cycle cannot be of length 3 (because $\lambda = 0$) or 4 (because $\mu=1$). If the cycle is of length 5, then the number of common vertices of $x$ and $y$ is $1$, otherwise it is $0$. Since $G$ is strongly regular by Lemma \ref{newlemma} part {\it 3} there exist non adjacent edges which belong to a $C_5$. So $L(G)$ is strongly regular if and only if any two non-adjacent edges belong to a $C_5$. \\

 The same reasoning applies if $\mu>1$, except that now $f_1$ and $f_2$ may belong to $C_4,C_5$ or $C_6$. If $f_1$ and $f_2$ belong to $C_4$ then the number of common vertices of $x$ and $y$ is $2$;  otherwise it is $1$ or $0$. Since $\mu>1$, in $G$ there are edges which belong to $C_4$. So $\Gamma(G)$ is strongly regular if and only if any two non-adjacent edges belong to a $C_4$.\\
\end{proof}

A similar result holds when $\lambda=1$.

\begin{thm}\label{20}
	If $G$ is a connected strongly regular graph with parameters $(n,k,1,\mu)$, then $\Gamma(G)$ is edge-regular. Furthermore, $\Gamma(G)$ is strongly regular if and only if $\mu>1, k=4,$ and any two non-adjacent edges belong to a $C_4$.
\end{thm}
\begin{proof}
	Since $G$ is strongly regular, $\Gamma(G)$ is $2(k-2)$-regular by Theorem \ref {15}. Now consider two adjacent vertices $v_1$ and $v_2$ in $\Gamma(G)$. If $e_1$ and $e_2$ be the corresponding edges in $G$, then the number of common vertices of $v_1$ and $v_2$ is same as the number of $K_{1,3}$'s in which $e_1$ and $e_2$ are present. By Lemma \ref{11}, since the number of such induced $K_{1,3}$ is $k-2$, any two adjacent vertices have $k-2$ common neighbours. Hence $\Gamma (G)$ is an edge-regular graph with parameters $(\frac{nk}{2},2k-4,k-2)$.\\

	To prove that $\Gamma(G)$ is strongly regular, consider two vertices $v_1$ and $v_2$ which are non-adjacent in $\Gamma(G)$. Let $e_1$ and $e_2$ be the corresponding edges in $G$. Then in $G$ either $e_1$ and $e_2$ span a triangle  or $e_1$ and $e_2$ are non-adjacent.\\

	In the first case, the number of common neighbours of $v_1$ and $v_2$ is same as the number of edges which form $K_{1,2}$ with both $e_1$ and $e_2$. By Lemma \ref{11} it is same as $k-2$.\\

	In the latter case, the number of common vertices is same as the number of induced $P_4$'s with end edges $e_1$ and $e_2$. To find the number of induced $P_4$'s we consider the following cases.\\
	1) $\mu=1$\\
	2) $\mu >1$\\
	If $\mu =1$ by Observation $3$, $e_1$ and $e_2$ may belong to a $C_5$ or $C_6$. By Observation $5$, $C_4$'s are forbidden in $G$. Therefore the number of such induced $P_4$ is one if $e_1$ and $e_2$ belong to a $C_5$; otherwise it is zero. Since $G$ contains $C_5$, $\Gamma(G)$ is strongly regular if and only if $k-2=1$ and any two non-adjacent edges belong to at least one $C_5$. But this is not possible since $k$ is an even number by Lemma \ref{11}.\\

	If $\mu>1$, by the same argument in the above theorem
	$\Gamma(G)$ is strongly regular if and only if $k-2=2$ and any two edges belong in a $C_4$. Hence $\Gamma(G)$ is strongly regular if and only if $k=4, \mu>1$ and any two non adjacent edges belong to a $C_4$.\\
\end{proof}

Furthermore, we can deduce a result on the connectivity of $\Gamma(G)$ under the same conditions as for the previous two results.

\begin{thm}\label{17'}
Let $G$ be a connected strongly regular graph with parameters $(n,k,\lambda,\mu)$. Suppose that $\lambda=0$ or $1$. Then $\Gamma(G)$ is $2$-connected.
\end{thm}
\begin{proof}
We have the following cases.\\

\ni{\bf{Case 1:}} $\lambda=0$\\
	Since $G$ is $K_3$-free, $\Gamma(G)\cong L(G)$. When we consider any two vertices of  $L(G)$, by Lemma \ref{newlemma} Part {\it 4}, the corresponding edges belong to a $C_4,C_5$ or $C_6$ (there can't be a vertex adjacent to all endpoints of the corresponding edges, since $G$ is $K_3$-free). So in $L(G)$ any two vertices belong to a  $C_4,C_5$ or $C_6$. Hence it is $2$-connected.\\

\ni{\bf{Case 2:}} $\lambda=1$\\
	In order to show that $\Gamma(G)$ is $2$-connected it is enough to show that any two vertices in $\Gamma(G)$ belong to a cycle. Consider two vertices $x$ and $x'$ in $\Gamma(G)$. We consider separately the following cases. \\
	(1) $xx'$ is an edge in $\Gamma(G)$.\\
	(2) $xx'$ is not an edge in $\Gamma(G)$.\\

\ni{\bf Subcase A:} $xx'$ is an edge in $\Gamma(G)$ implies that the corresponding edges $e$ and $e'$ are incident in $G$ and do not belong to a $K_3$. By Lemma \ref{11} there exists an edge $l$ such that $e$ and $l$ span a $K_3$ in $G$. Similarly there exists edge $l'$ such that $e'$ and $l'$ span a $K_3$ in $G$. Then the vertices corresponding to $e,l',l,e'$ form a path $P_4$ in $\Gamma(G)$. Then this path together with the edge $xx'$ form a $C_4$ in $\Gamma(G)$.\\

Figure \ref{Gallai_cycle} shows the edges used in this argument to form a cycle in $\Gamma(G)$. The reader may wish to take the time to thoroughly digest this argument, as variants of it will be used throughout the next subcase. \\

\begin{figure}
\begin{centering}
\includegraphics[width=7cm,height=6cm,keepaspectratio]{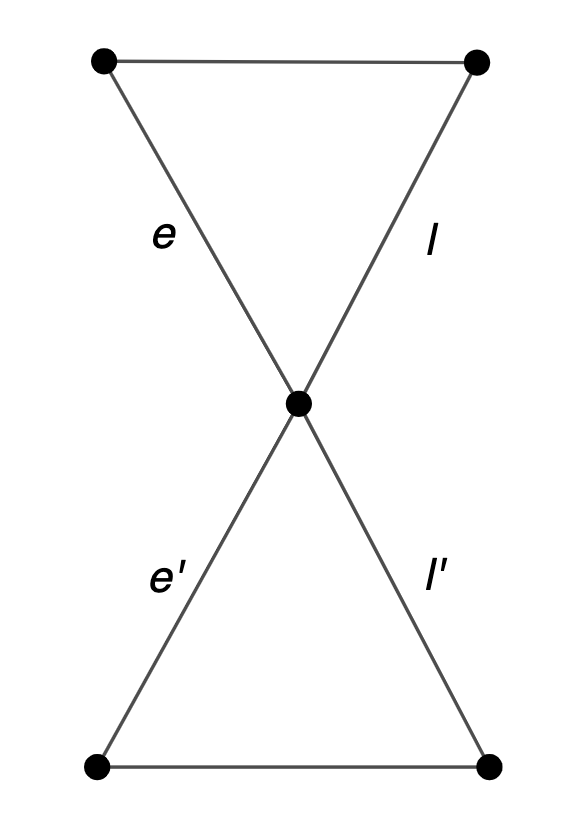}
\par\end{centering}
\caption{Edges forming a $C_4$ in $\Gamma(G)$} \label{Gallai_cycle}
\end{figure}

\ni{\bf{Subcase B:}} $xx'$ not an edge in $\Gamma(G)$, which means that either the corresponding edges are incident in $G$ and  belong to a $K_3$ or are not incident. For the first case, by Lemma \ref{11}, there exist at-least two edges $l$ and $l'$ which span a $K_3$ in $G$. Then the vertices corresponding to $e,l,e',l',e$ span a $C_4$ in $\Gamma(G)$. For the latter case by Lemma \ref{newlemma} part {\it 4}, either $e$ and $e'$ belong to a cycle $C_n$ where $3<n<7$ in $G$, or there is a vertex adjacent to all endpoints of $e,e'$. Let us consider the cycle case first. If it is an induced $C_n$ then clearly $x$ and $x'$ belong to a $C_n$ in $\Gamma(G)$. If not, there are edges which belong to a $K_3$ in $G$. Consider two edges $e_1$ and $e_2$ with  a common vertex $u$ and span a $K_3$ in $G$. By Lemma \ref{11} there exist another edges $l$ and $l'$ incident on $u$ which span a $K_3$ in $G$. Since diamonds are forbidden in $G$, the vertices corresponding to $e_1,l, e_2$ form a $P_3$ in $G$. So, if two edges span a $K_3$ in $G$ then by the above explanation we can find an edge  such that the vertices corresponding to the edges of $K_3$ and  the new edge form a $P_3$ in $G$. Also if $e_1,e_2,e_3$ are three consecutive edges in $C_n$ and if $e_1$ and $e_2$ span a $K_3$ in $G$, since diamonds are forbidden in $G$, $e_2$ and $e_3$ cannot span a $K_3$ in $G$. So by the above explanation in $\Gamma(G)$, we can find a cycle $C_n$ of length at most $9$  containing the vertices $x$ and $x'$. That is, in $\Gamma(G)$ any two vertices belong to at least one  $C_n$. Hence $\Gamma(G)$ is $2$-connected.\\

Finally, if there is a vertex $v$ adjacent to each endpoint of $e$ and $e'$ then a similar argument allows us to build a cycle in $\Gamma(G)$, as follows. Since $k \geq 4$ each endpoint of $e$ is adjacent to a point other than $v$, and this can't be an endpoint of $e'$ since diamonds are forbidden. This new edge cannot span a triangle with $e$ or containing $v$, since $\lambda = 1$. Similar new edges can be found at the endpoints of $e'$. The edges $ef_1s_1s'_2f'_2e'f'_1s'_1s_2f_2e$ form a $C_{10}$ in $\Gamma(G)$.
\end{proof}

As discussed in the introduction, the spectra of Gallai graphs is of considerable interest, and one of the most important properties of strongly regular graphs is simple form of their spectra, and the ease with which it can be calculated. This is the object of the following theorem, which is well-known.

\begin{thm}\cite{Bap}\cite{cv2}\label{13}
	The distinct eigenvalues of a connected strongly regular graph $G$ with parameters $(n,k,\lambda,\mu)$ are $k,\dfrac{1}{2}[(\lambda-\mu)+\sqrt{(\lambda-\mu)^2+4(k-\mu)}]$ and $k,\dfrac{1}{2}[(\lambda-\mu)+\sqrt{(\lambda-\mu)^2-4(k-\mu)}]$.\\
\end{thm}

We will now work through a series of examples that illustrate our results, and will in each case indicate the spectra of the Gallai graphs in question.

\begin{example}\label{17} \normalfont
	Let $G$ be a connected strongly regular $C_6$-free graph with parameters $(n,k,0,1)$. Since $\mu=1$, by Lemma \ref{newlemma} parts {\it 4} and {\it 5} any two non-adjacent edges must belong to either $C_6$ or $C_5$ or both (the other possibility in Lemma \ref{newlemma} part {\it 4}, a vertex adjacent to all endpoints of both edges, cannot occur since $\lambda=1$). Since the graph is $C_6$-free any two edges must belong to $C_5$, hence by Theorem \ref{16}, $\Gamma(G)$ is a strongly regular graph with parameters $(\dfrac{nk}{2},2k-2,k-2,1)$. By Theorem \ref{13}, the distinct eigenvalues of $\Gamma(G)$ are $ 2k-2,\dfrac{1}{2}[(k-3)+\sqrt{(k-3)^2+4(2k-3)}]$ and $\dfrac{1}{2}[(k-3)-\sqrt{(k-3)^2+4(2k-3)}]$.
\end{example}

\begin{example} \label{17a} \normalfont
	Let $G$ be a connected strongly regular graph with parameters $(n,k,0,\mu)$, where $\mu>1$. Suppose also that every pair of non-adjacent edges belong to a $C_4$. Then, since $\mu>1$, by Theorem \ref{16}, $\Gamma(G)$ is a strongly regular graph with parameters $(\dfrac{nk}{2},2k-2,k-2,2)$.The distinct eigenvalues of  $\Gamma(G)$ are therefore $2k-2,\dfrac{1}{2}[(k-4)+\sqrt{(k-4)^2+4(2k-4)}]$ and $\dfrac{1}{2}[(k-4)-\sqrt{(k-4)^2+4(2k-4)}]$ by Theorem \ref{13}.
\end{example}

\begin{example}\label{21} \normalfont
	Let $G$ be a connected strongly regular graph with parameters $(n,k,1,\mu)$, where $\mu>1$ and any two edges are contained in a $C_4$. by Theorem \ref{20}, $\Gamma(G)$ is a connected strongly regular graph with parameters $(2n,4,2,2)$, and then by Theorem \ref{13} the distinct eigenvalues of  $\Gamma(G)$ are $4,\sqrt{2}$ and $-\sqrt{2}$.
\end{example}

\subsection{Anti-Gallai graphs}

The following is the analogue of Theorem \ref{14} for anti-Gallai graphs.

\begin{thm} \label{31}
		Let $G$ be a connected graph. The anti-Gallai graph $\Delta(G)$ is disconnected if and only if there exists a partition of the edge set of $G$ into $E_1, E_2, ... E_p$ where $p\geq2$, such that if $e_i\in E_i$ and $e_j\in E_j$ are incident in $G$ then $e_i$ and $e_j$ do not belong to a $K_3$, and there exists at least one pair of this type.
\end{thm}
\begin{proof}
Suppose that $\Delta(G)$ is disconnected and let $\Delta_1, \Delta_2,...\Delta_p$ with $p\geq 2$ be the components of $\Delta(G)$. As in the proof of Theorem \ref{14} consider

$$E_i=\{ e\in G : e \mbox{ is an edge corresponding to a vertex $v$ in $\Delta_i$}\},$$

where $1\leq i \leq p$. Clearly, $E_i $ is a partition for  $E(G)$. Since the  connectedness of $G$ implies the  connectedness of $\L(G)$,  at least one edge $e_i \in E_i$ is incident with some $e_j\in E_j$ with $j \neq i$. However, $\Delta_i$ and $\Delta_j$ are different components of $\Delta(G)$, and hence if $e_i\in E_i$ is incident with $e_j\in E_j$ then they do not belong to a triangle in $G$.\\

For the converse assume that such a partition exists for $E(G)$. Then for any $i$ and $j$, the vertices corresponding to the edges in $E_i$ and $E_j$ induce different components in $\Delta(G)$.\\
\end{proof}

We can deduce from this a result on strongly regular graphs with $\lambda = 2$, as follows.
\begin{thm} \label{33}
	Let $G$ be a connected strongly regular graph with paramenters $(n,k,2,\mu)$, where each vertex in $G$ belongs to exactly one wheel. Then $\Delta(G)$ is connected and edge-regular.
\end{thm}
\begin{proof}
Suppose $\Delta(G)$ is disconnected. Then by Theorem \ref{31} there exists at least two edges  $e_i\in E_i$ and $e_j\in E_j$ which are incident but do not belong to a $K_3$ in $G$. Let $u$ be the common vertex of both $e_i$ and $e_j$. By assumption, since the neighbouring vertices induce a wheel in $G$, there exist edges $e_{i+1},e_{i+2},..., e_{j-1},e_j$ such that the pairs of edges $(e_i,e_{i+1}), (e_{i+1},e_{i+2}),...,(e_{j-1},e_j)$ belongs to a $K_3$ in $G$. Since $e_i\in E_i$, the edges $e_{i+1},e_{i+2},..., e_{j-1},e_j$ all belong to $E_i$. This is a contradiction, and hence $\Delta(G)$ is connected.\\

\ni By Theorem \ref{15}, $\Delta (G)$ is a $4$-regular graph. Since $G$ is not $K_4$ any two adjacent vertices in $\Delta (G)$ have only one common neighbour. Hence $\Delta(G)$ is edge-regular with parameters $(\frac{nk}{2},4,1)$.\\
\end{proof}

As in the previous subsection, we will now work through a series of examples illustrating our results.

	\begin{example} \normalfont
		Let $G$ be a connected strongly regular graph with parameters $(n,k,0,\mu)$. Since $\lambda=0$, $G$ is $K_3$-free. Hence $\Delta(G)$ is totally disconnected. The spectrum of $\Delta(G)$ is $(0^{\frac{kn}{2}})$.
	\end{example}

\begin{example} \normalfont
	Let $G$ be a connected strongly regular graph with parameters $(n,k,1,\mu)$. Since $\lambda=1$, every edge of $G$ belong to exactly one $K_3$ and no two $K_3$ share a common edge. Therefore $\Delta(G)$ is the disjoint union of $\dfrac{kn}{6}$ triangles. $\Delta(G)$ therefore has the spectrum $(-1^{\frac{kn}{3}}, 2^\frac{kn}{6})$.
\end{example}

\begin{example} \normalfont
	 Let us consider the anti-Gallai graph of the line graph of the complete bipartite graph, $L(K_{n,n})$. $L(K_{n,n})$ contains $2n$ copies of $K_n$ sharing common vertices, where two $K_n$'s have no common edges and the edges of two copies of $K_n$ do not belong to a $K_3$. Hence $\Delta(L(K_{n,n}))$ is the disjoint union of $2n$ copies of $\Delta (K_n)$. Since every pair of edges in $K_n$ spans a triangle, $\Delta(K_n)\cong L(K_n)$, and $L(K_n)$ is a well-known family of graphs known as the {\it triangular graphs}. These are strongly regular with parameters $(\frac{n(n-1)}{2},2(n-2), n-2, 4)$ and spectrum $(2(n-2)^1,(n-4)^{n-1}, -2^{n(n-3)/2})$. Hence the spectrum of $\Delta(L(K_{n,n}))$ is $(2(n-2)^{2n},(n-4)^{2n(n-1)}, -2^{n^2(n-3)})$.
\end{example}

Our final result applies to a situation in which we are not able to identify $\Delta(G)$ easily but in which we can say something about the spectrum only.

\begin{defn}
Suppose $\alpha_1 \leq \ldots \leq \alpha_n$ and $\beta_1 \leq \ldots \leq \beta_m$, with $m < n$. We will say the sequence of $\beta$'s {\it interlace} the sequence of $\alpha$'s if $\alpha_k \leq \beta_k \leq \alpha_{k+n-m}$ for $k = 1, \ldots, m$.
\end{defn}

The following two lemmas are famous results in spectral graph theory.

	\begin{lem}\cite{cv}\label{35}
		If $G$ be a  graph and $H$ an induced subgraph. Then the eigenvalues of $H$ interlace those of $G$.
	\end{lem}

\begin{lem} \cite{Bal} \label{36}
	If $G$ is a $k$-regular graph, then $spec(G)\in[-k,k]$, and $k \in spec(G)$.
\end{lem}

We need two further definitions.

\begin{defn}
 Let $G_1$ and $G_2$ be vertex-disjoint graphs. The \textit{join} of $G_1$ and $G_2$, denoted $G_1\vee G_2$, is the supergraph of $G_1+G_2$ in which each vertex of $G_1$ is adjacent to every vertex of $G_2$ \cite{Bal}.
\end{defn}

\begin{defn}
The {\it semi-total point graph} $R(G)$ of a graph $G$ is obtained from $G$ by adding a new vertex corresponding to every edge of $G$, then joining each new vertex to the end vertices of the corresponding edge i.e; each edge of $G$ is replaced by a triangle \cite{Sam}.
\end{defn}

\begin{figure}
\begin{centering}
\includegraphics[width=7cm,height=6cm,keepaspectratio]{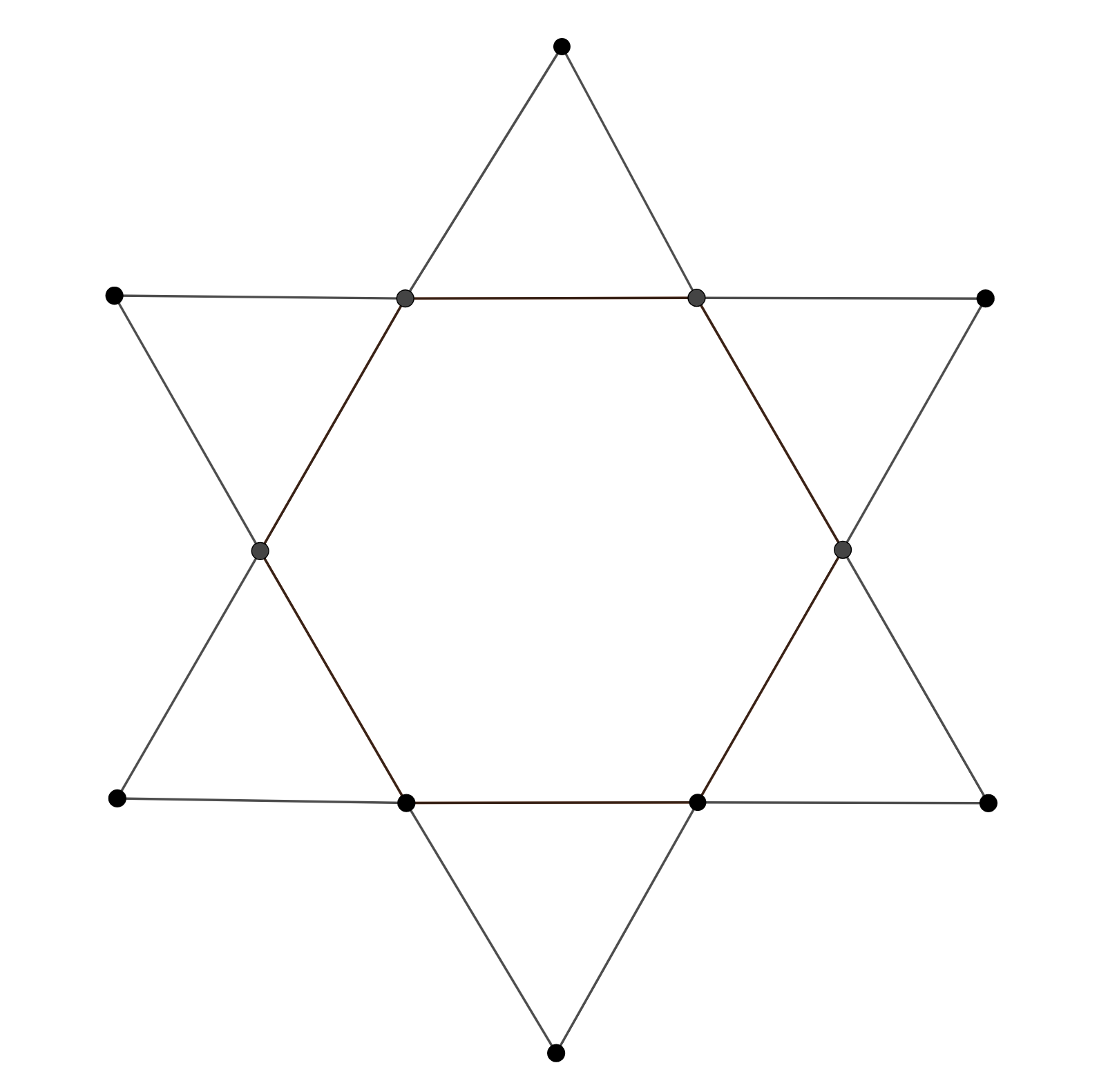}
\par\end{centering}
\caption{$R(C_6)$} \label{semitot}
\end{figure}

The semi-total point graph of a cycle, $R(C_n)$, will be of special interest to us, and Figure \ref{semitot} represents this graph for $n=6$.
The following lemma connects joins, semi-total point graphs, and anti-Gallai graphs.

\begin{lem}\cite{Jee2}\label{37}
If $G= H\vee K_1$, where $H$ is $K_3$-free, then $\Delta(G)$ is the semi-total point graph of $H$.
\end{lem}

We are now prepared to prove the following theorem.

\begin{thm}
	Let $G$ be a connected strongly regular graph with parameters $(n,k,2,\mu)$, and suppose that each vertex in $G$ belongs to exactly one wheel. Then  $spec(\Delta(G)) \subseteq [-4,4]$ and the eigenvalues of $R(C_k)$ interlace those of $\Delta(G)$.
\end{thm}	
\begin{proof}
By Theorem \ref{33}, $\Delta(G)$ is a $4$-regular graph. Therefore, by Lemma \ref{36} the largest eigenvalue of $\Delta(G)$ is $4$ and $spec(G)\in [-4,4]$. Also, by assumption each vertex belong to exactly one wheel ($C_n \vee K_1$).  Since $\lambda=2$  a vertex in $G$ together with its neighbours form an induced $k$-wheel. Then, by Lemma \ref{37}, it is clear that $\Delta(G)$ contains the semi-total point graph  $R(C_k)$ as an induced subgraph.  Therefore, by Lemma \ref{35}, the eigenvalues of $R(C_k)$ interlace those of $\Delta(G)$.\\
\end{proof}

To apply this result requires knowledge of the spectrum of $R(C_k)$. The eigenvectors of this graph can be found in a similar manner to that of a cycle, namely by placing powers of a complex $n$-th at the points of the cycle, and then finding what the remaining values must be at the other points. We omit the details, but the reader may check that that spectrum of $R(C_k)$ is the values of the form $\frac{r \pm \sqrt{r^2+ 4r + 8}}{2}$, where $r=2\cos(2\pi k/n)$ for $k=1,\ldots, n$.

\section{Acknowledgements}

We would like to thank Misha Lavrov for helpful conversations.



\begin{thebibliography}{}

\bibitem{PHRC} Anand, P., Escaudro, H., Gera, R., Martell, C., Triangular line graph and word sense disambiguation, Discrete Applied Mathematics 159 (2011), 1160 - 1165.

\bibitem{Apa1} Aparna Lakshmanan S., Characterization of some special classes of the Gallai and the anti-Gallai graphs, Discourse 1 (2013), 85 - 89.

\bibitem{Apa} Aparna Lakshmanan S., Rao, S. B., Vijayakumar, A., Gallai and anti-Gallai graphs of a graph, Math. Bohem., 132(1) (2007), 43 - 54.

\bibitem{Bal}  Balakrishnan, R., Ranganathan, K., A text book of graph theory, Springer (1999).

\bibitem{Bap} Bapat, R.B., Graphs and Matrices, Springer.

\bibitem{Bro}  Brouwer, A., Van Maldeghem, H., Strongly regular graphs, https://homepages.cwi.nl/\~aeb/math/srg/rk3/srgw.pdf.

\bibitem{cam} Cameron, P., Strongly regular graphs, Topics in Algebraic Graph Theory, 102 (2004), 203-221.

\bibitem{cv} Cvetcovic, D., Graphs and Their Spectra, (Thesis), Univ.Beograd,Publ. Elektrotehn. Fak. Ser. Mat. Fiz.. 354-356 (1971), 1 - 50.

\bibitem{cv1} Cvetcovic, D., Doob, M., and Sachs, H., Spectra of Graphs - Theory and application, academic press, New York, 1980.

\bibitem{cv2} Cvetcovic, D., Peter Rowlinson and Slobodan Simic, An Introduction to the Theory of Graph Spectra, London Mathematical Society

\bibitem{gal} Gallai, T., Transitiv orientierbare graphen, Acta Math. Hung., 18(1-2) (1967), 25--66.

\bibitem{garg} Garg, P., Sinha, D., Goyal, S., Eulerian and hamiltonian properties of Gallai and anti-Gallai total graphs, J. Indon. Math. Soc. 21 (2) (2015), 105 - 116.

\bibitem{godroy} Godsil, C., Royle, G., Algebraic graph theory. Vol. 207. Springer Science and Business Media, 2001.

\bibitem{Jee1} Palathingal, J.J., Aparna Lakshmanan S., Gallai and anti-Gallai Graph Operators, Electronic Notes in Discrete Mathematics 6.3 (2017), 447 - 453.

\bibitem{Jee3} Palathingal, J.J., Indulal, G., Aparna Lakshmanan S., Spectrum of Gallai Graph of Some Graphs, Indian J. of Pure and Appl. Math. 51(4)(2020), 1829 - 1841.

\bibitem{Jee2} Palathingal, J.J., Indulal, G., Aparna Lakshmanan S., Spectrum of anti-Gallai Graph of Some Graphs, Indian J. of Pure and Appl. Math., (to appear).

\bibitem{IG} Gutman, I., Suriha, I., On the Nullity of Line Graph of Tree, Discrete Math. 232 (2001), 35 - 45.

\bibitem{joos} Joos, F., Le, V.B., Rautenbach, D., Forests and trees among Gallai graphs, Discrete Math. 338(2)  (2015), 190 - 195.

\bibitem{Ka} Kavitha, K., David, N.G., Dominator coloring of some classes of graphs, International Journal of Mathematical Archive, 3(11) (2012), 3954 - 3957

\bibitem{Pri} Prisner, E., "Graph Dynamics", Longman, 1995.

\bibitem{Le1} Le, V. B., Gallai graphs and Anti-Gallai Graphs, Discrete Math., 159 (1996), 179 - 189.

\bibitem{Le2} Le, V. B., Mortality of Iterated Gallai Graphs, Period. Math. Hungar., 27(2) (1993), 105 - 124.

\bibitem{L}  Beineke, L., Wilson, R., Cameron, P.,Topics in Algebraic Graph Theory: 102 (Encyclopedia of Mathematics and its Applications, Series Number 102), Cambridge University Press, 2004.

\bibitem{agnes1} Poovathingal, A., Kureethara, J., Deepthy, D., A survey of the studies on Gallai and anti-Gallai graphs, Communications in Combinatorics and Optimization, 6 (1) (2021), 93-112.


\bibitem{agnes2} Poovathingal, A., Kureethara, J., Some characterizations of {G}allai graphs, AIP Conference Proceedings, 2236 (1) (2020), 060003.

\bibitem{Sam} Sampathkumar, E., Chikkodimath, S.B., The Semitotal graphs of a graph-II, J. Karnatak Univ. Sci, 18 (1973), 281 - 284.


\end{thebibliography}
\end{document}